\documentclass[11pt]{amsart}
\usepackage{amscd,amssymb,amsmath,enumerate}
\usepackage{amsmath,amsthm}
\usepackage[all]{xy}
\usepackage{graphics}
\usepackage[pdftex]{graphicx}
\usepackage{latexsym}
\usepackage[all]{xy}
\usepackage{psfrag}
\xyoption{matrix} \xyoption{arrow}
\usepackage{parskip}
\usepackage{hyperref} 

\newtheorem{proposition}{Proposition}[section]
\newtheorem{theorem}[proposition]{Theorem}

\newtheorem{definition}[proposition]{{Definition}}
\newenvironment{defn}{\begin{definition} \rm}{\end{definition}}

\newtheorem{remark}[proposition]{{Remark}}

\newtheorem{Example}[proposition]{Example}

\newcommand{\cQ}{{Q}}

\newcommand{\cS}{{\mathcal S}}

\newcommand{\Hom}{\operatorname{Hom}\nolimits}

\newcommand{\rad}{\operatorname{rad}\nolimits}

\usepackage{color}
\definecolor{candyapplered}{rgb}{1.0, 0.03, 0.0}
\makeatletter
\def\thm@space@setup{%
  \thm@preskip=0.7cM \thm@postskip=0.3cM
}
\makeatother

\begin{document}

\date{today}
\title[A survey]
{Special multiserial algebras, Brauer configuration algebras and more : a survey}

\author[Green]{Edward L.\ Green}
\address{Edward L.\ Green, Department of
Mathematics\\ Virginia Tech\\ Blacksburg, VA 24061\\
USA}
\email{green@math.vt.edu}

\author[Schroll]{Sibylle Schroll}
\address{Sibylle Schroll\\
Department of Mathematics \\
University of Leicester \\
University Road  \\
Leicester LE1 7RH, UK
}
\email{schroll@le.ac.uk}

\subjclass[2010]
{16S37, 14M99,16W60}

\keywords{}
\thanks{This work was supported through the UK Engineering and Physical Sciences Research Council grant\href{http://gow.epsrc.ac.uk/NGBOViewGrant.aspx?GrantRef=EP/K026364/1}{EP/K026364/1}. }

\begin{abstract} We survey results on multiserial algebras, special multiserial algebras and  Brauer configuration algebras.   A structural property of modules
over a  special multiserial algebra is presented.  Almost gentle algebras are introduced and we  describe some
results related to this class of algebras.   We also report on the structure of radical cubed zero symmetric
algebras. 
\end{abstract}
\date{\today}
\maketitle

\section{Introduction}\label{sec-intro}
In this paper we present a brief survey of the results of  \cite{GS, GS3, GS4, GS2}. These
papers represent an attempt to generalize biserial algebras, special biserial algebras, and Brauer graph algebras.  We show that the generalizations satisfy properties that 
are analogous to the properties of the algebras they generalize. Biserial, special
biserial, and Brauer graph algebras have been extensively studied,
see, for example,  \cite{A1,AAC,BS,CBVF,Du,E,Fuller,G,GreenSchrollSnashall,
GreenSchrollSnashallTaillefer,K,MarshSchroll,Roggenkamp,Schroll, SW,SkowronskiSurvey,WW}  and also the survey articles \cite{Schroll2,Schroer} and the references therein.   

The generalizations of these algebras are the multiserial algebras  \cite{KY, VHW}, 
the special multiserial algebras \cite{GS2, VHW} and
the Brauer configuration algebras \cite{GS}.   These algebras are typically wild,  but as our
results indicate, they seem to
be  worthwhile algebras to study.     For  proofs of the results
presented in this survey, examples, and further details, the
reader is referred to \cite{GS, GS3, GS4, GS2}.

After presenting the definitions of multiserial algebras and special multiserial algebras,
we define Brauer configuration algebras and algebras defined by cycles.  We then
look at how these algebras are related  to each other in Theorems \ref{thm-multi}, \ref{thm-conn}, and
\ref{thm-quot}.

Like biserial algebras, gentle algebras have received a great deal of attention,  see for example 
\cite{ABCP,BS,CC, GR,Schroll}.  In \cite{GS3} we define  almost gentle algebras, a natural generalization of gentle algebras,
and show that trivial extensions of such algebras by the dual of the algebra
are Brauer configuration algebras, Theorem \ref{thm-triv}.

The most surprising result concerns modules over a  special multiserial algebra.  We define multserial modules and show that, even though  special multiserial algebras
tend to be wild, every finitely generated module over a  special multiserial algebra
is a multiserial module, Theorem \ref{multiserial theorem}.

There have recently been substantial advances in the classification of tame symmetric radical cubed zero algebras \cite{Benson, ErdmannSolberg1, SkowronskiSurvey}.  We show that over an algebraically
closed field, every indecomposable symmetric radical cubed zero algebra is a Brauer configuration algebra,
Theorem \ref{thm-r30}.

Throughout this paper, $K$  denotes a field, $\cQ$ a finite
quiver with vertex set $\cQ_0$, and arrow
set $\cQ_1$, and $K\cQ$ the path algebra.  We let $J$ be the ideal in
$K\cQ$ generated by the arrows of $\cQ$.   We assume, for simplicity, that
every algebra occurring in the survey is of the form $K\cQ/I$ with
$J^N\subseteq I\subseteq J^2$, for some $N\ge 3$ and  $\rad^2 (K\cQ/I) \neq 0$.
We also assume
that the algebras are indecomposable as rings.

\section{Multiserial and special multiserial algebras}\label{sec-multi}  We begin with a brief
discussion of biserial algebras and special biserial algebras.  Recall that a $K$-algebra $\Lambda$ is \emph{biserial} if for every indecomposable projective
left or right module $P$, there are uniserial left or right modules
$U$ and $V$, such that $\rad(P) = U+ V$ and $U \cap V $ is either zero or simple.  

A $K$-algebra $\Lambda$ is a \emph{special biserial} algebra if the following two
properties hold:
\begin{enumerate}
\item[S1] For each arrow $a\in\cQ_1$, there is at most one arrow $b\in\cQ_1$ such that
$ab\notin I$, and there is at most one arrow $c\in\cQ_1$ such that $ca\notin I$.
\item[S2] For each vertex $v\in\cQ_0$, there are at most 2 arrows starting at $v$ and
at most 2 arrows ending at $v$.

\end{enumerate}

A third type of algebra,  that is closely related, is a Brauer graph algebra which is constructed from a Brauer graph.
A \emph{Brauer graph} is a 4-tuple $\Gamma=( \Gamma_0,\Gamma_1,\mu,\mathfrak o)$ where
$(\Gamma_0,\Gamma_1)$ is a finite  (undirected) graph where loops and multiple edges are permitted, and
$\mu\colon \Gamma_0\to \mathbb Z_{\ge 1}$  is a function, called the
\emph{multiplicity function}.  To describe $\mathfrak o$, called the \emph{orientation}
of $\Gamma$, we need a few more definitions.
A vertex $v\in\Gamma_0$ is a \emph{truncated vertex} if $v$ is an endpoint of exactly one edge  and
$\mu(v)=1$.  All other vertices are called \emph{nontruncated vertices}.  The orientation
$\mathfrak o$ is a choice, for each nontruncated vertex $v$, of a cyclic ordering  of the
edges having $v$ as an endpoint.   Since the construction of a Brauer graph algebra is
a special case of the construction of a Brauer configuration algebra, we postpone the
construction until later in the paper.  It should be noted that Brauer graphs and Brauer graph
algebras generalize Brauer trees and Brauer tree algebras, which appear in the study
of group rings of finite  and tame representation type \cite{ Janusz, Dade, Donovan}.
We have the following connections.

\begin{theorem}\label{thm-biser}\cite{ Roggenkamp, Schroll, SW} Let $\Lambda=K\cQ/I$ with  $\rad^2(\Lambda) \neq 0$.
If $\Lambda$ is a special biserial algebra, then $\Lambda$ is biserial.
Furthermore, if $K$ is algebraically closed 	then 
$\Lambda$ is a symmetric special biserial algebra if and only if $\Lambda$ is a Brauer graph
algebra.  
\end{theorem}

One of the most important features of special biserial algebras  is that their representation theory has a combinatorial framework in terms of strings and bands where the finitely  generated indecomposable modules are given by  string and band modules \cite{WW}.  Thus
 special  biserial algebras are tame \cite{WW}.  The algebras we now define are usually of wild representation
type.

Let $\Lambda$ be a $K$-algebra and let $M$ be a finitely generated $\Lambda$-module. We say that $M$ is {\it left multiserial} (respectively {\it right multiserial}) if there is a positive integer $n$ and left (respectively right) uniserial modules
$U_1,\dots,U_n$, such that $\rad(M) = \sum_{i=1}^nU_i$ and, if $i\ne j, U_i \cap U_j $ is either zero or simple. 
We say that  $\Lambda$ is \emph{multiserial} if
 every indecomposable projective
 left or right module is left or right multiserial \cite{GS}, see also  \cite{VHW, KY}.  

 We say $\Lambda=K\cQ/I$ is \emph{special multiserial}
 if condition (S1) holds, that is if:
\begin{enumerate}
\item[]
  For each arrow $a\in\cQ_1$, there is at most one arrow $b\in\cQ_1$ such that
$ab\notin I$, and there is at most one arrow $c\in\cQ_1$ such that $ca\notin I$.
\end{enumerate}

Since this is one of the defining properties of a special biserial algebra, a special biserial
algebra is a special multiserial algebra.
The next result generalizes part of Theorem \ref{thm-biser}.

\begin{theorem}\label{thm-multi}\cite{GS} A special multserial algebra is a multiserial
algebra.
\end{theorem}

\section{Algebras defined by cycles}\label{sec-cycles}

Before turning our attention to Brauer configurations  and Brauer configuration algebras, we introduce another new class of algebras, {\it  algebras defined by cycles \cite{GS2}}.

 Let $\cQ$ be a quiver. A \emph{simple cycle} at a vertex $v$ is a cycle  $C=a_1a_2\cdots a_n$ in $\cQ$  where
the $a_i$ are arrows, $a_1$
starts at $v$, $a_n$ ends at $v$ and there are no repeated arrows in $C$,  that is $a_i \neq a_j$ for all $i \neq j$.   The cyclic permutations
of $C$ are $a_1\dots a_n$,  $a_2a_3\cdots a_na_1, \dots, a_na_1a_2\cdots a_{n-1}$.
A \emph{defining set of cycles} is a set $\cS$ of simple cycles of length at least 1, such that
\begin{enumerate}

\item[D1] if $ C\in\cS$, then every cyclic permutation of $C$ is in $\cS$,
\item[D2] every arrows occurs in some $C\in\cS$, and
\item[D3] if an arrow occurs in two cycles in $\cS$, then   the cycles are a cyclic permutation
of each other.

\end{enumerate} 
Note  that some of the cycles can be  loops.
If $\cS$ is a defining set of cycles, then a \emph{multiplicity function for $\cS$} is a function
 $\nu \colon\cS\to \mathbf Z_{\ge 1}$ such that if $C,C'\in\cS$ and $C$ is a cyclic permutation of $C'$,
then $\nu (C)=\nu (C')$.  Furthermore, we require  that,   if $C\in\cS$ is a loop, then $\nu (C)>1$.  We call $(\cS,\nu )$ a \emph{defining pair}.

If $(\cS,\nu )$ is a defining pair, then define $I_{(\cS,\nu )}$ to be the ideal in $K\cQ$ generated  by all elements of the following 3 types:
\begin{enumerate}
\item[Type  1] $C^{\nu (C)}-C'^{\nu (C')}$, if $C,C'$ are cycles in $\cS$  at some
vertex $v$,
\item[Type 2] $C^{\nu (C)}a_1$, if $C=a_1a_2\cdots a_m$,
\item[Type 3] $ab$, where $a,b\in\cQ_1$ and  $ab$ is not a subpath of $C^2$, for all
$C\in\cS$  such that $C$ is not a loop.

\end{enumerate}

\begin{defn}
Let $Q$ be a quiver, $(\cS,\nu)$  a defining pair on $Q$ and $I_{(\cS,\nu )}$ the ideal defined by $(\cS,\nu)$.  Then
we call $\Lambda=K\cQ/I_{(\cS,\nu)}$  the \emph{algebra defined by  $(\cS,\nu)$} and we
say $\Lambda$ is 
\emph{defined by cycles}.
\end{defn}

\section{Brauer configurations and Brauer configuration algebras}\label{sec-bca}

We now turn to Brauer configurations.   A Brauer configuration is a 4-tuple of the
form $\Gamma=(\Gamma_0,\Gamma_1,\mu,\mathfrak o)$.  The elements of
the finite set $\Gamma_0$ are called \emph{vertices} and the elements of
$\Gamma_1$ are called \emph{polygons}.  If $V\in\Gamma_1$, then $V$ is
a finite multiset of vertices, that is, the elements of $V$ are vertices, possibly
with repetitions.  The function $\mu\colon\Gamma_0\to \mathbf Z_{\ge 1}$ is
called the\emph{ multiplicity function}.   A vertex $\alpha\in\Gamma_0$ is called \emph{truncated}
if $\mu(\alpha)=1$, $\alpha$ occurs  exactly once in one polygon $V$  where $V=\{\alpha,\beta\}$ and
\[ \text{(the number of total occurences of }\beta\text{ in polygons counting repetitions)}\cdot(\mu(\beta))\ge 2  .\]
If a vertex is not truncated it is \emph{nontruncated}.
 The \emph{orientation} $\mathfrak o$ is a choice, for each nontruncated vertex $\alpha\in\Gamma_0$, of a cyclic ordering of the polygons that contain $\alpha$, counting
repetitions.  

In \cite{GS} we called the above ``reduced'' Brauer configurations. There are some technical restrictions imposed by the assumption that the associated Brauer configuration
algebra, (to be defined), is indecomposable.  These restrictions can be found in \cite{GS}.

An example will help to clarify the above definition.  Suppose $\Gamma =(\Gamma_0,\Gamma_1,\mu,\mathfrak o)$ where
$\Gamma_0=\{1,2,3,4\}, \Gamma_1=\{V_1,V_2,V_3\} $ where $V_1=\{1,2,3\},
V_2=\{1,1,3\}$, $V_3=\{2,3\}$ and $V_4 = \{2,4\}$.  Let $\mu(i)=1$, for $i=1,\dots,4$.  The only truncated
vertex is 4.  The orientation $\mathfrak o$ is a choice for each nontruncated vertex of an ordering
of the polygons that contain the vertex.
For example, for 1, the ordering might be $V_1<V_2<V_3<V_2(<V_1)$ with the last inequality
because it is a cyclic ordering. Any cyclic permutation, for example, $V_2<V_3<V_2<V_1(<V_2)$
is considered the same ordering.  On the other hand,  the ordering at 1 could also have been chosen to be
 $V_1<V_2<V_2<V_3(<V_1)$ which yields a different configuration.   

A Brauer configuration is a generalization of a Brauer graph since every Brauer graph is
a Brauer configuration with the restriction that
every  polygon is a set  with 2 vertices.

There are ``realizations" of 
a Brauer configuration  obtained by attaching actual polygons to the elements in $\Gamma_1$
with the vertices of the polygon labelled by the appropriate elements of $\Gamma_0$.  Then one
identifies vertices of the actual polygons that have the same label.  Another way to realize a Brauer configuration is as a hypergraph where the set of vertices is given by $\Gamma_0$  and where $\Gamma_1$ corresponds to the set of hyperedges \cite{GS3}.

We now define the Brauer configuration algebra $\Lambda_{\Gamma}$ associated to a Brauer configuration
$\Gamma=(\Gamma_0,\Gamma_1,\mu,\mathfrak o)$.  The quiver $\cQ_{\Gamma}$ of
$\Lambda_{\Gamma}$ has vertex set $\Gamma_1$.  If $\Gamma_1=\{V_1,\dots, V_n\}$,
we write $(\cQ_{\Gamma})_0=\{v_1,\dots,v_n\}$ to distinguish between the polygons in $\Gamma$
and the vertices in $\cQ_{\Gamma}$.   For every nontruncated vertex $\alpha\in\Gamma_0$, there is
a cyclically ordered sequence of polygons containing $\alpha$ given by the orientation $\mathfrak o$.  For each nontruncated $\alpha
\in\Gamma_0$,  we view the ordered sequence of polygons as a  simple cycle $C_{\alpha}$ in $\cQ_{\Gamma}$.
If $\alpha,\beta\in\Gamma_0$ are nontruncated with $\alpha\ne \beta$, then we require that
 $C_{\alpha}$ and $C_{\beta}$ 
have no arrows in common.  Note that if $\alpha$ is in exactly one polygon $V$ and $\mu(\alpha)>1$, then $C_{\alpha}$ is a loop at $v$ in $\cQ_{\Gamma}$.

The arrow set of $\cQ_{\Gamma}$ are exactly 
the arrows occuring in the $C_{\alpha}$, for $\alpha$ a
nontruncated vertex in $\Gamma$.  Set $\cS=\{C_{\alpha}\mid \alpha \text{ a nontruncated
vertex in }\Gamma\}$.  Define  $\nu\colon\cS\to\mathbb Z_{\ge 1}$ by $\nu(C_{\alpha})
=\mu(\alpha)$,  where, on the  right hand side, $\mu$ is the multiplicity function of $\Gamma$. 
It is easy to see that  $(\cS,\nu)$ is a defining pair.

The Brauer configuration algebra $\Lambda_{\Gamma}$ associated to $\Gamma$ is
the algebra  $K\cQ_{\Gamma}/I_{(\cS,\nu)}$ defined by $(\cS,\nu)$.
This description of $\Lambda_{\Gamma}$ is slightly different than the original
description given in \cite{GS} but both  yield isomorphic algebras.

\section{Connection results}

The following results connect all the algebras.

\begin{theorem}\label{thm-conn} \cite{GS4, GS2}  Let $K$ be an algebraically closed field  and let $\Lambda=K\cQ/I$ with $\rad(\Lambda)^2\ne 0$, be
a symmetric $K$-algebra. Then the following statements are equivalent.
\begin{enumerate}
\item $\Lambda$ is a symmetric special multiserial algebra.
\item $\Lambda$ is defined by  cycles.
\item $\Lambda$ is a Brauer configuration algebra.
\end{enumerate}
\end{theorem}

Another result generalizing a result about  special  biserial algebras \cite{WW} is the following.

\begin{theorem}\label{thm-quot}\cite{GS2}  Let $K$ be a field.  Every special multiserial $K$-algebra is  a quotient of
a symmetric special multiserial $K$-algebra.
\end{theorem}


\section{Examples}\label{sec-ex}

\begin{Example}\label{ex1.pdf}{\rm Let $\cQ$ be the quiver

\begin{center}\includegraphics[scale=.5]{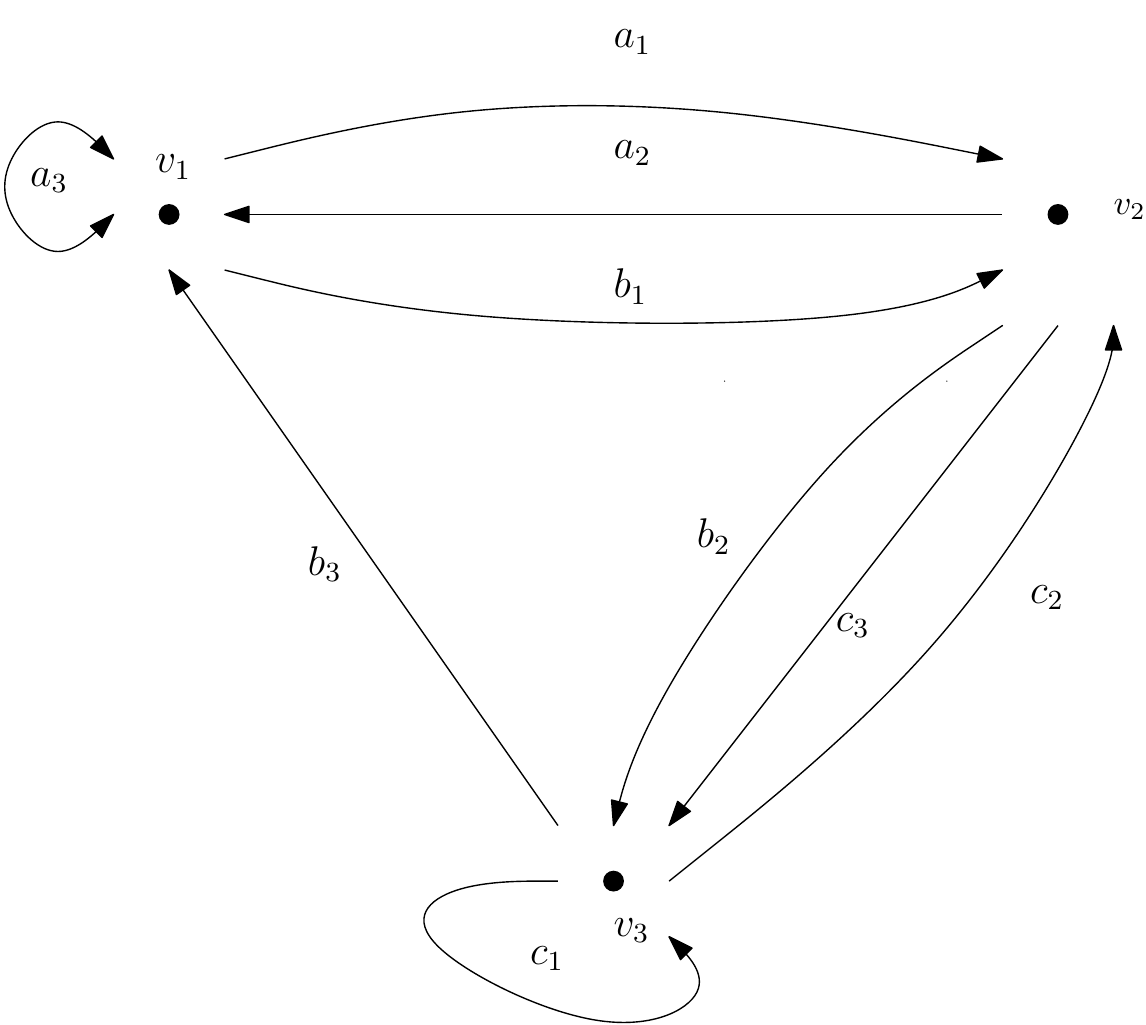}\end{center}
Let $\cS$ consist of the following cycles and their cyclic    permutations.:
$a_1a_2a_3, b_1b_2b_3, c_1c_2c_3$ and define $\nu\colon\cS\to\mathbb  Z_{\ge 1}$
by $\nu\equiv 1$.  We see that $(\cS,\nu)$ is a defining pair. Let $\Lambda$ be defined
by $(\cS,\nu)$.  In particular, a generating set for the ideal of relations for  $\Lambda$
is given in Section \ref{sec-cycles}.  By Theorem \ref{thm-conn}, $\Lambda $ is a Brauer configuration algebra and a special
multiserial algebra.  The Brauer configuration is $\Gamma=(\Gamma_0, \Gamma_1,\mu,\mathfrak o)$ where $\Gamma_0=\{1,2,3\}$, $\Gamma_1=\{V_1,V_2,V_3\}$ with $V_1=\{1,1,2\},
V_2=\{2,3,3\}$, and $V_3=\{1,2,3\}$.  The multiplication function $\mu$  is identically 1 and
the orientation is given below. 

A realization of the Brauer configuration is
\begin{center}\includegraphics[scale=.4]{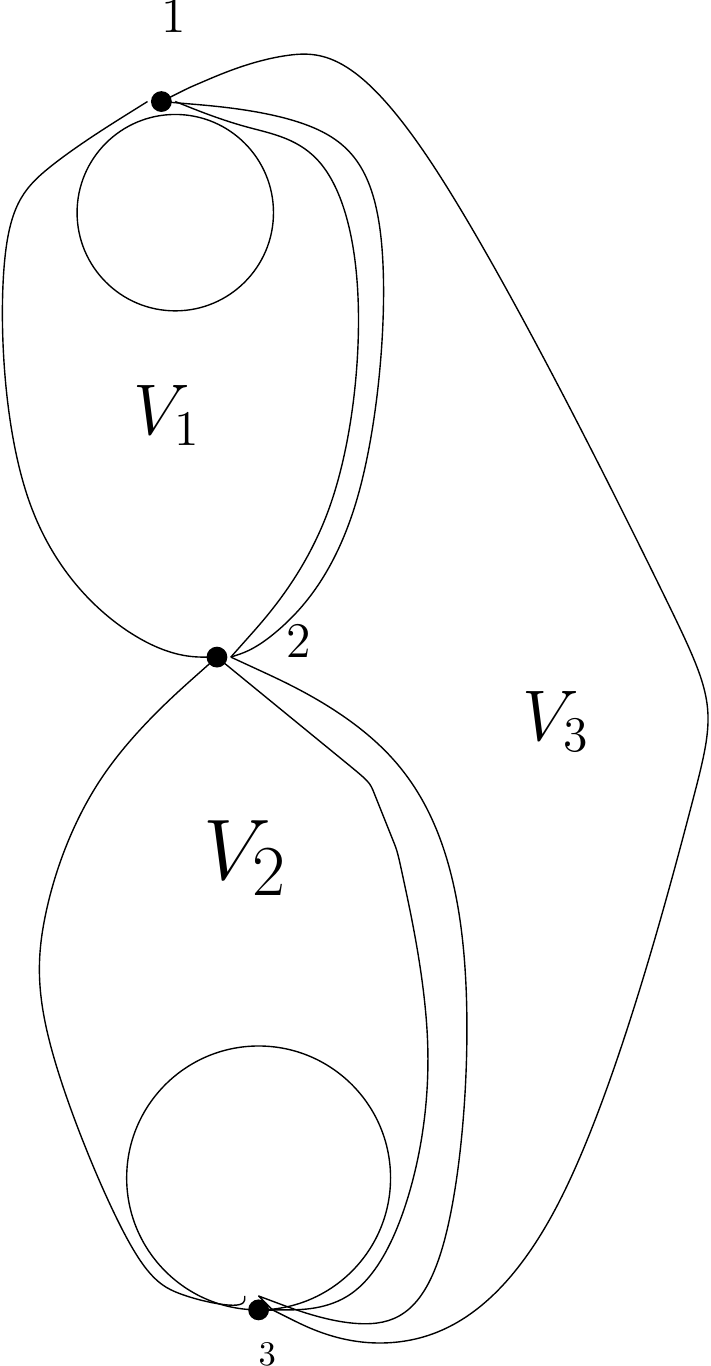}
\end{center}with clockwise orientation at each vertex.
}\end{Example}

\begin{Example}\label{ex-ex2}{\rm
Let the quiver be
\begin{center}\includegraphics[scale=.8]{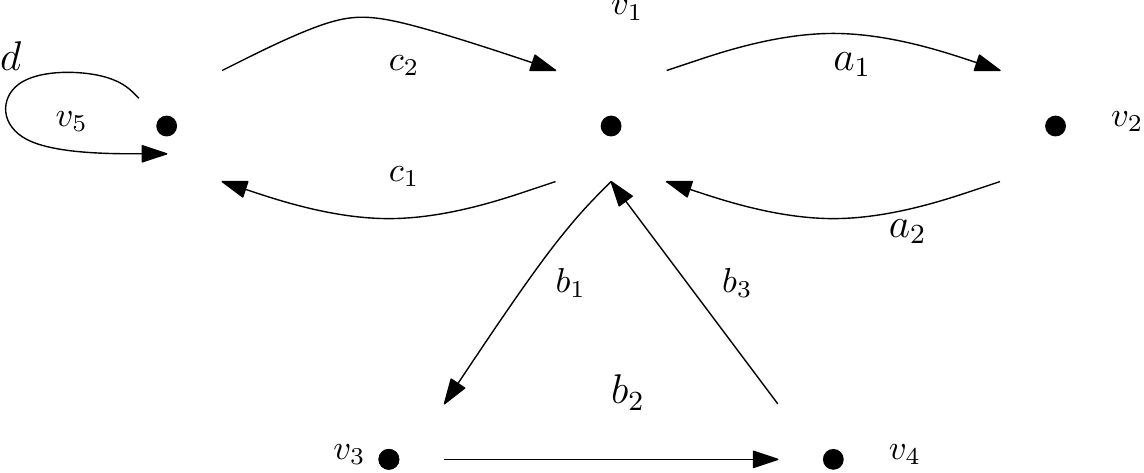}\end{center}
Let $\cS$  be the set of simple cycles consisting of all the
permutations of the cycles $a_1a_2, b_1b_2b_3, c_1c_2,d$  with $\nu : \cS \to \mathbb{Z}_{\geq 1}$ defined by
$\nu(a_1a_2)=1=\nu( b_1b_2b_3)=\nu( c_1c_2)$, and $\nu(d)=2$. We see that
$(\cS, \nu)$ is a defining  pair. 
Take $\Lambda$ to be the algebra defined by $(\cS,\nu)$.

One can check that $\Lambda$ is isomorphic to the Brauer configuration
algebra associated to the Brauer configuration $\Gamma=(\Gamma_0,\Gamma_1,\mu,
\mathfrak 0)$ with $\Gamma_0=\{1,2,3,4,5,6,7\}$, $\Gamma_1=\{V_1,V_2,V_3,V_4,V_5\}$
with $V_1=\{1,2,3\}, V_2=\{1,4\}, V_3=\{2,5\}, V_4=\{2,6\}, V_5=\{3,7\}$.  The multiplicity
function is 1 for all vertices except for vertex 7 where $\mu(7) =2$.  We see that
vertices 4,5,6 are the truncated vertices.
The orientation is given by the clockwise  ordering of the polygons at each nontruncated
vertex in the realization of the Brauer configuration  below.

A realization of the Brauer configuration is:
\begin{center}\includegraphics[scale=.8]{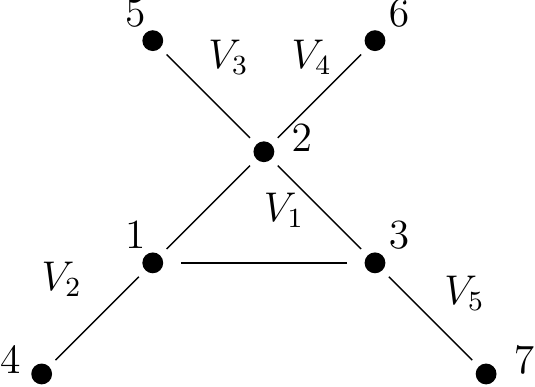}\end{center}
}\end{Example}
\section{Almost gentle algebras}
Recall that an algebra $\Lambda=K\cQ/I$ is \emph{gentle} if it is special biserial, that is if (S1) and (S2) hold, and if additionally   S2$^{\rm op}$ holds, that is for any arrow $a$ there exists at most one arrow $b$ such that $ab \in I$ and at most one arrow $c$ such that $ca \in I$ and if  $I$ 
is  generated by paths of length 2.  

 An obvious way to generalize the notion of a gentle algebras is:  An algebra
$\Lambda=K\cQ/I$ is \emph{almost gentle} if it is special multiserial and if $ I$ can be
generated by paths of length 2.  

 A characterisation of gentle algebras is through the trivial extension by their minimal injective co-generator. That is, an algebra $A$ is gentle if and only if the trivial extension $T(A) = A \ltimes \Hom_k(A,k) $ is special biserial \cite{PS}.

We have the following results.
\begin{theorem}\label{thm-triv}  \cite{GS3}  If $A$ is an almost gentle algebra then $T(A)$ is
a symmetric special multiserial algebra  and hence it is isomorphic to a Brauer configuration algebra. If $T(A)$ is isomorphic to the Brauer configuration
algebra associated to the Brauer configuration  $\Gamma=(\Gamma_0,\Gamma_1,\mu,\mathfrak o)$,
then $\mu\equiv 1$.

\end{theorem}

There is a converse. Let $\Lambda=K\cQ/I$ be a symmetric special multiserial algebra.
Then by Theorem \ref{thm-conn}, there is a defining pair $(\cS,\nu)$ such that
$\Lambda$ is defined by $(\cS,\nu)$. Assume that $\nu\equiv 1$. We say a set of arrows $D$ is an \emph{admissible
cut} if $D$ consists of one arrow for each permutation class of cycles in $\cS$.  
Let $\widehat \cQ$ be the quiver obtained from $\cQ$ by removing the arrows in $D$.
Let $\widehat I$ be the ideal $I\cap K\widehat\cQ$.

\begin{theorem}\label{thm-adm}\cite{GS3} Let $\Lambda=K\cQ/I$ be a symmetric special multiserial algebra  viewed as an algebra defined by cycles with defining pair $(\cS,\nu)$ where $\nu\equiv 1$.
Keeping the notation of the previous paragraph,  let $D$ be an admissible cut in $\cQ$.
Then $K\widehat\cQ/\widehat I$ is an almost gentle algebra. 
\end{theorem}

If $A$ is an almost gentle algebra, then there is an admissible cut in $T(A)$ such that
the construction above yields an algebra isomorphic to $A$.  Moreover,
if $\Lambda$ is a symmetric special  multiserial algebra (with $\nu\equiv 1$) and $D$ an admissible cut,
then 	$T(K\widehat\cQ/\widehat I)$ is isomorphic to $\Lambda$.  The proofs of these results are 
found 
in \cite{GS3}.

\section{Representations of multiserial algebras}
As mentioned earlier,  special multiserial algebras are typically wild.  Thus, it is surprising that
there is a structural result about all finitely generated modules over such an algebra.
 Recall that a left (respectively, right) finitely generated module is multiserial if it is left multiserial (respectively, right multiserial), that is if $\rad(M)$ as a left (resp., right) module is a sum of uniserial modules that intersect pairwise in the zero module or in a simple module. 

\begin{theorem}\label{multiserial theorem}
 Let 
$A$ be a special multiserial $K$-algebra, and $M$  an indecomposable finitely generated
left (respectively, right)  $A$-module. Then $M$ is multiserial. 
\end{theorem} 

\section{Radical cubed  zero}\label{sec-r30}

In this section we present some results about symmetric algebras with the condition that
  their Jacobson radical cubed is  zero.  The proof of the next result requires an algebraically closed field.

\begin{theorem}\label{thm-r30} \cite{GS4}
Let $K$ be an algebraically closed field 
and let $A=K\cQ/I$ be a  finite dimensional  basic indecomposable $K$-algebra. Suppose that $A$ is symmetric and  that
$\rad^3(A)=0$ but $\rad^2(A)\ne 0$.
Then $A$ is isomorphic to a Brauer configuration algebra.
\end{theorem}

It is well-known that there are equivalence relations on 
 symmetric matrices with non-negative integer
entries,  as well as on finite undirected graphs, such that there is a one-to-one correspondence between
the respective equivalence classes, see \cite{GS} for details.  To each equivalence
class  one can associate a radical cubed zero symmetric algebra.  Now consider the  radical cubed zero Brauer configuration algebras  with Brauer configuration $\Gamma=
(\Gamma_0,\Gamma_1,\mu,\mathfrak o)$ that satisfy 
the following  2 properties:
\begin{enumerate}
\item there are no repeated vertices occurring in any  polygon in $\Gamma_1$,
\item if $\alpha\in \Gamma_0$ then either $\alpha$ is in exactly one polygon and
$\mu(\alpha)=2$ or $\alpha$ is in exactly two polygons and $\mu(\alpha)=1$.

\end{enumerate}
 We show that there is an equivalence relation on radical cubed zero Brauer configuration algebras with  properties (1) and (2) such that the
equivalence classes are in one-to-one
correspondence with the radical cubed zero symmetric algebras associated to the equivalence
classes of finite graphs \cite{GS}.

\end{document}